\documentstyle{article}
   \font\tenmsb=msbm10
   \font\sevenmsb=msbm7
   \font\fivemsb=msbm5

\newfam\msbfam
      \textfont\msbfam=\tenmsb
      \scriptfont\msbfam=\sevenmsb
      \scriptscriptfont\msbfam=\fivemsb
\def\Bbb#1{{\fam\msbfam #1}}

\newcommand\qed{{\hspace*{\fill}Q.E.D.\vskip12pt plus 1pt}}

\newcommand\sE{{\cal E}}
\newcommand\sJ{{\cal J}}
\newcommand\sL{{\cal L}}
\newcommand\sO{{\cal O}}

\newcommand\gra{\alpha}
\newcommand\grg{\gamma}
\newcommand\grl{\lambda}

\newcommand\reals{{\Bbb R}}
\newcommand\comp{{\Bbb C}}
\newcommand\pn[1]{{\Bbb P}^{#1}}
\newcommand\opn[1]{\sO_{{\Bbb P}^{#1}}}
\newcommand\proof{\noindent{\em Proof.}\ \ }

\newtheorem{theorem}{Theorem}[section]
\newtheorem{lemma}[theorem]{Lemma}
\newtheorem{proposition}[theorem]{Proposition}

\begin{document} 

\title{Line bundles for which a projectivized jet bundle is a product}  
\author{Sandra Di Rocco and Andrew J. Sommese}
\date{June 25, 1998}
\maketitle

\begin{abstract} 
We characterize the triples $(X,L,H)$, consisting of line bundles $L$ 
and $H$ on a complex projective manifold $X$, such that for some  
positive integer $k$, the $k$-th holomorphic jet bundle of $L$, $J_k(X,L)$, 
is isomorphic to a direct 
sum $H\oplus\dots\oplus H$.

\noindent {\em 1991 Mathematics Subject Classification}. 14J40, 14M99.\newline
\indent{\em Keywords and phrases.}   jet bundle, complex projective manifold,
projective space, Abelian variety.
\end{abstract}

\section*{Introduction}
Let $X$ be a complex projective  manifold. A large amount of   
 information on the geometry of an embedding $i: X\hookrightarrow  \pn N$  
 is contained in the ``bundles of jets" of the line bundles $k\sL=i^*\sO_{\pn{N}}(k)$ 
for $k\geq 1$. 
The $k$-th jet bundle of a line bundle $L$ (sometimes called 
the $k$-th principal part of $L$) is usually denoted by $J_k(X,L)$, or by
$J_k(L)$ when the space $X$ is clear from the context.
The bundle $J_k(k\sL)$ is spanned by the $k$-jets of global sections of $k\sL$. If $\gra:\pn{}(J_1(\sL))\to \pn{N}$ is the map given by the $1$-jets of elements of
$H^0(X,\sL)$, then $\gra(\pn{}(J_1(\sL))_x)={\Bbb T}_x$, where ${\Bbb T}_x$ 
is the embedded tangent space at $x$. Similarly $\pn{}(J_k(k\sL))_x$ is 
mapped to the ``$k$-th embedded tangent space" at $x$ (see \cite{GrHa} for more details).
Given this interpretation of the projectivized jet bundle, it is natural to expect that 
a  projectivized jet bundle is a product of the base and 
a projective space only under very rare circumstances. 

In this paper we analyze the more general setting where $L$ is a line bundle on $X$ (with no hypothesis on its positivity) and $J_k(L)=\oplus H$, with  $H$
 a line bundle on $X$.   Some years ago the second author \cite{So} analyzed the
pairs $(X,L)$ under the stronger hypothesis that $J_k(L)$ is a trivial bundle. Though, 
the results in this paper do not follow from \cite{So}, the only possible projective 
manifolds  with a projectivized jet bundle of some line bundle being  a 
product of the base manifold and a projective space,  
turn out to be the same, i.e., Abelian varieties 
and projective space. We   completely characterize
 the possible triples $(X,L,H)$.  In fact, on Abelian varieties we have the necessary 
and sufficient condition that
$L\in {\mbox{\rm Pic}_0}(X)$, and if $L=\sO_{\pn{n}}(a)$, then only the range 
$a\geq k$ or $a\leq -1$ can occur. 

A significant amount of the research on this paper was done during our 
stay at Max-Planck-Institut f\"ur Mathematik in Bonn, to which we would like to express our gratitude for the 
excellent working conditions.  The second author would also like to thank the 
Alexander von Humboldt Foundation for its support.

\section{Some preliminaries}
We follow the usual notation of algebraic geometry.
 Often we denote the direct sum of  $m$ copies of a vector bundle $\sE$, for some integer $m> 0$,
by $\displaystyle\bigoplus_m\sE$. We freely use the additive notation for line bundles.  We often use the same symbols for a vector bundle and its associated
locally free sheaf of germs of holomorphic sections. We say a vector bundle
is spanned if the global sections generate each fiber of the bundle.
For general references  we refer to \cite{Book} and \cite{KS}.

By $X$ we will always denote a projective manifold over the complex numbers $\comp$ 
and by $L$ a holomorphic line bundle on $X$.

Let $k$ be a nonnegative integer. The {\em $k$-th jet bundle} $J_k(L)$ associated to $L$ is defined as
the vector bundle of rank ${k+n}\choose{n}$ associated to the sheaf 
$p^*L/(p^*L\otimes\sJ_{\Delta}^{k+1})$, where $p:X\times X\to X$ is
the projection on the first factor and $\sJ_{\Delta}$ is the sheaf of ideals of the diagonal, $\Delta$, of $X\times X$.
Let $j_k:H^0(X,L)\times X \to J_k(L)$ be the map which associates to each section $s\in H^0(X,L)$ 
its $k$-th jet. That means that locally, choosing coordinates $(x_1,...,x_n)$ and a 
trivialization of $L$ in a neighborhood of $x$, $j_k(s,x)=(a_1,...,a_{{n+k}\choose {n}})$,
 where the $a_i$'s are the coefficients of the terms of degree up to $k$ in the Taylor
 expansion of $s$ around $x$. Notice that the map $j_k$ is surjective if and only if $H^0(X,L)$
 generates all the $k$-jets at  all points $x\in X$. For example $j_1$ being 
surjective is equivalent to $|L|$ defining an immersion of $X$ in $\Bbb P^{h^0(L)-1}$.

We will often use the associated exact sequence:
$$
0\to T_X^{*(k)}\otimes L\to J_k(L)\to J_{k-1}(L)\to 0\;\;\;\;\;\;\;\;(j_k)
$$
where $T_X^{*(k)}$ denotes the $k$-th symmetric power of the cotangent bundle of $X$. 
There is an injective bundle map \cite[p.\ 52]{KS}
$$\grg_{i,j}: J_{i+j}(L)\to J_i(J_j(L)).$$

Using the sequence $(j_k)$ it is easy to see that
$$\det J_k(L)=\frac{1}{n+1}{{n+k}\choose{n}}(kK_X+(n+1)L)$$

\begin{lemma}\label{split}
Let $X$ be a compact K\"ahler variety and $L$ a holomorphic line 
bundle on $X$. Then $c_1(L)=0$ in $H^1(T^*_X)$ 
 if and only if the bundle sequence $(j_1)$ splits. If $c_1(L)= 0$ in $H^1(T^*_X)$,
then $J_k(L)\cong J_k(\sO_X)\otimes L$.
\end{lemma}
\proof A local computation, using Cech coverings, shows that the Atiyah 
class defined by the sequence $(j_1)$ is the cocycle $c_1(L)\in H^1(X,T_X^*)$.  Thus
$c_1(L)=0$ in $H^1(T^*_X)$  if and only if the bundle sequence $(j_1)$ splits.

If   $c_1(L)=0$ in $H^1(T^*_X)$ then $L$ has constant transition functions.
 From this, the isomorphism $J_k(L)\cong J_k(\sO_X)\otimes L$ follows.
\qed

Let $A$ and $B$ be vector bundles on $X$, we recall that if $\gra\in H^1(A\otimes B^*)$ represents the vector bundle extension $E$, then any nonzero multiple $\grl\gra$ gives an isomorphic extension. If $\grl=0$ this is of course false since we would get the trivial extension.
It follows that:
\begin{lemma}\label{ext}
Let $\grl$ be a nonzero integer, then $J_1(\grl L)\cong J_1(L)\otimes (\grl -1)L$
\end{lemma}
\proof Consider the extension :
$$0\to T_X^*\otimes \grl L\to J_1(L)\otimes (\grl -1)L\to\grl L\to 0\;\;\;\;\;\;\;\;(j_1)\otimes (\grl -1)L$$
represented by $c_1(L)+ c_1((\grl -1)L)=\grl c_1(L)=c_1(\grl L)$.
Then, from what was observed above, $J_1(\grl L)$, which is the vector 
bundle extension given by $c_1(\grl L)$, is isomorphic to $ J_1(L)\otimes (\grl -1)L$.
\qed

\section{The basic examples}
In this section we will characterize line bundles with splitting $k$-th 
jet bundles on  Abelian varieties and on $\pn{n}$. These  turn out to be 
 the only possible examples.

 \begin{proposition}\label{Abelian} Let $L$ and $H$ be line bundles on an Abelian variety $X$, then the following assertions are equivalent:
\begin{itemize}
\item $J_k(L)\cong \oplus H$ for some $k$;
\item $H\cong L$ and $L\in {\mbox{\rm Pic}_0}(X)$.
\end{itemize}
In particular $J_k(L)$ splits in the sum of spanned line bundles only when
$L=\sO_X$ and $J_k(L)=\oplus \sO_X$.
\end{proposition}
\proof Let $X$ be an Abelian variety and assume that $J_k(L)=\oplus H$ for some line bundle $H$. Then the sequences $(j_m)$ for $m\le k$ imply that there is a surjection of the trivial bundle $J_k(L)\otimes (-H)$ onto $L-H$.
Thus the line bundle $L-H$ is a spanned line bundle with $${{n+k}\choose{n}}c_1(L-H)=0,$$
which implies $L=H$. Using this we can find a direct summand $L$ of $J_k(L)$ whose
image in $J_1(L)$ maps onto $L$ under the map $J_1(L)\to L\to 0$
of the sequence $(j_1)$.  Thus the sequence $(j_1)$ splits 
 and therefore $c_1(L)=0$ in $H^1(T^*_X)$ by Lemma \ref{split}.

Conversely assume that we have a line bundle $L\in {\mbox{\rm Pic}_0}(X)$.
Then by Lemma \ref{split}, we have that $J_k(L)\cong J_k(\sO_X)\otimes L$.
Using the triviality of $T_X^{(j)}$ for all $j\ge 0$, it follows that $J_k(\sO_X)$ is trivial. \qed

\begin{proposition}\label{pn} $J_k(\opn{n}(a))$ is isomorphic to a direct sum 
$\displaystyle \bigoplus_{{k+n}\choose{n}}\opn{n}(q)$ for some $k,a,q\in {\Bbb Z}$
with $k>0$ if and only if $q=a-k$ and either  $a\geq k$ or $a\leq -1$. 
\end{proposition}
\proof Let $L=\opn{n}(a)$. First notice that if the $k$-th jet bundle  splits as $\displaystyle J_k(L)=\bigoplus_{{k+n}\choose{n}}\opn{n}(q)$ for some nonnegative $k$ and some integer $q$, then 
$$\det J_k(L)=\opn{n}\left({{n+k}\choose{k}}(a-k)\right),$$ which implies that $q=a-k$.
 Then the fact that the $k$-th jet bundle of $\opn{n}(a)$ always has sections for $a\ge 0$
 rules out the cases $a=0,...,k-1$.

Lemma \ref{ext} gives $J_1(\opn{n}(a))\cong J_1(\opn{n}(1))\otimes \opn{n}(a-1) 
$. Then $J_1(\opn{n}(1))=\oplus \opn{}$ (see, e.g., \cite{So}) implies
$$J_1(\opn{n}(a))=\bigoplus_{n+1}\opn{n}(a-1)$$
Dualizing $\grg_{1,1}$ and using $$J_1\left(J_1(\opn{n}(a)\right)=J_1\left(\bigoplus_{n+1}\opn{n}(a-1)\right)=\bigoplus_{n+1}\left(\bigoplus_{n+1}\opn{n}(a-2)\right)$$
 we obtain the quotient
$$\bigoplus_{n+1}\left(\bigoplus_{n+1}\opn{n}(a-2)\right)^*\to J_2(\opn{n}(a))^*\to 0$$
and thus, tensoring by $\opn{n}(a-2)$
$$\bigoplus_{(n+1)^2}\opn{n}\to J_2(\opn{n}(a))^*\otimes \opn{n}(a-2)\to 0$$
Then the vector bundle $J_2(\opn{n}(a))^*\otimes \opn{n}(a-2)$ is spanned with
$\det(J_2(\opn{n}(a))^*\otimes \opn{n}(a-2))=\opn{n}$
which implies that  $J_2(\opn{n}(a))^*\otimes \opn{n}(a-2)=\oplus\opn{n}$.
Iterating this argument one gets  $\displaystyle  J_k(\opn{n}(a))^*\otimes \opn{n}(a-k)=\bigoplus_{{k+n}\choose{n}}\opn{n}.$ \qed

\section{The main result}
In this section we characterize complex projective manifolds having a line bundle whose projectivized $k$-jet bundle is a product of the base manifold 
and a projective space. 
\begin{theorem}
Let $L$ and $H$ be holomorphic line bundles on $X$. Then $J_k(L)=H\oplus\dots\oplus H$ if and only if the triple $(X,L,H)$ is one of the two below:
\begin{enumerate}
\item  $(X,L,L)$ where $X$ is Abelian and $L\in {\mbox{\rm Pic}_0}(X)$
\item  $(\pn{n}, \opn{n}(a),\opn{n}(a-k))$ and $a\geq k$ or $a\leq -1$.
\end{enumerate}
\end{theorem}
\proof Assume $J_k(L)=H\oplus\dots\oplus H$. Tensoring the sequence $(j_k)$ by $H^*$ gives the quotient
$ \displaystyle \bigoplus_{{n+k}\choose{n}}\sO_X\to L\otimes H^*\to 0.$
Tensoring the sequence $(j_k)$ by $H^*$  and then dualizing it
gives the quotient
$ \displaystyle \bigoplus_{{n+k}\choose{n}}\sO_X\to T_X^{(k)}\otimes (H\otimes L^*)\to 0.$
It follows that $L\otimes H^*, T_X^{(k)}\otimes (H\otimes L^*)$ and thus $T_X^{(k)}$ are spanned bundles over $X$. 

First assume that  the canonical bundle $K_X$ is nef.
Since $T_X^{(k)}$ is spanned, we conclude that $\det T_X^{(k)}$, which
is a spanned negative multiple of $K_X$
is trivial. Thus    $T_X^{(k)}$ is a trivial bundle. From this we conclude that $L-H$ is trivial.  It follows that the trivial bundle $J_k(L)^*\otimes L$ has a filtration with quotient bundles $T^{(j)}_X$ for $0\le j\le k$.  This shows that $T^{(j)}_X$ is trivial for all $j> 0$. Thus under the assumption
that $K_X$ is nef, we conclude that $X$ would be an Abelian variety.  Proposition \ref{Abelian} implies also that $L\in {\mbox{\rm Pic}_0}(X)$.
This gives the first case of the theorem.

 We are thus left with the case when $K_X$ is not nef. 
 Since 
$T_X^{(k)}$ is spanned $-K_X$ is nef. 
The cone theorem then yields the existence of an extremal ray $\reals_+[ \grg]$, with $1\leq -K_X\cdot\grg\leq n+1.$ Let $l$ be the normalization of $\grg$ and let
$$T_{X|l}=\bigoplus \sO_l(a_i),\;\;\;(L\otimes H^*)_l=\sO_l(b),$$
where by abuse of notation we denote the pullback of a bundle $\sE$ on $\gamma$ to $l$
by $\sE_l$.
Writing for simplicity  $T^{(k)}_{X|l}=(\oplus_i \sO_l(ka_i))\oplus P$ we get
$$T^{(k)}_{X|l}\otimes (L\otimes H^*)^*_l=(\oplus_i \sO(ka_i-b))\oplus (P\otimes\sO_l(-b))$$
Since $T^{(k)}_{X|l}\otimes (L\otimes H^*)^*_l$ is spanned we conclude that $ka_i\geq b\ge 0$. Note that $b> 0$. Indeed if $b=0$, then
$$0=\deg \det J_k(L)_l\otimes (-H_l)= \frac{1}{n+1}{{n+k}\choose{n}}kK_X\cdot l\not=0.$$
Thus   $a_i>0$ for all $i$'s. Moreover from sheaf injection
$\displaystyle 0\to T_l\to T_{X|l}$
we see that $T_{X|l}$ must contain a factor $\sO_l(a_i)$ with $a_i\geq 2$.
Then $\displaystyle -K_X\cdot\grg=\sum_{i=1}^{n+1}a_i\leq n+1$ 
implies that $a_j=2$  for one $j$ and $a_i=1$ for $i\neq j$, i.e.,  $-K_X\cdot\grg=n+1$.
Now from Mori's proof of Hartshorne conjecture (see \cite[\S 4]{L}) we  
  deduce that $X$ must be $\pn{n}$. At this point we have recovered case (2) by applying Proposition \ref{pn}.

Propositions \ref{Abelian} and \ref{pn}  show that if we are in cases (1) and (2) respectively, then $\displaystyle J_k(L)=H\oplus\dots\oplus H$ satisfied. \qed

{
\begin{tabular}{ll} Sandra Di Rocco&Andrew J. Sommese\\
 Department of Mathematics&Department of Mathematics\\
KTH, Royal Institute of Technology\ \ \ &University of Notre Dame\\
 S-100 44 Stockholm, Sweden&Notre Dame, Indiana 46556, U.S.A.\\
e-mail: sandra@math.kth.se&e-mail: sommese@nd.edu\\
http://www.math.kth.se/$\sim$sandra&http://www.nd.edu/$\sim$sommese
\end{tabular}
}

\end{document}